\baselineskip=14pt

\bigskip
{\bf STABLE EXTENSIONS BY LINE BUNDLES}
\bigskip

\bigskip

Let C be a generic curve of genus g. Denote by $U(n,d)$ the moduli space 
of stable vector bundles of rank n and degree d on C. Let $E$ be an element
of  $U(n,d)$. Consider the set of subbundles $E'$ of rank $n'$ of  $E$.
Let $E''$ be the corresponding quotient and $n''$ its rank.
 We have an exact sequence
$$(0,1)0\rightarrow E'\rightarrow E\rightarrow E''\rightarrow 0$$
 Define the integer
$$(0,2)s_{n'}(E)=min\{ n'deg(E'')-n''deg(E')\}=n'd-nmax\{deg(E')\}$$
Where the minimum or maximum is taken as $E'$ varies in the set of
 subbundles of 
$E$ of rank $n'$ . Then, $s_{n'}(E)\equiv n'd (n)$.
 It is known (cf.[N],[M-S]) that $0\le s_1(E)\le n''n'g$.
Consider now the following stratification of $U(n,d)$: denote by 
$$U_{s,n'}(n,d)=\{ E \in U(n,d)|s_{n'}(E)\le s\} , 0\le s\le n'n''g$$
As the function $s_{n'}(E)$ is upper semicontinuos in a family of
vector bundles, the loci $U_s(n,d)$ are closed in $U(n,d)$.
The expected dimension of $U_s(n,d)$ is given by (see [L] section 4)
$$d_{s,n'}=(n^2-n'n'')(g-1)+s+1$$ 
Non-emptiness of $U_{s,n'}(n,d)$ will imply in many cases that the
expected dimension is in fact the dimension (see [L] p.455).
Lange conjectured (cf.[L]) that the loci $U_{s,n'}$ should be non-empty
for $0\le s \le n'n''(g-1)$
In this paper prove the conjecture in the case $n'=1$ and we shall often drop
the subindex 1.

\bigskip 
\proclaim (0,3) Theorem. On a generic curve of genus g 
$dimU_{s,1}(n,d)=U_s(n,d)=d_{s,1}, 0< s\le(n-1)g, s\equiv d(n)$.
\bigskip

We are going to prove the result by a degeneration argument. 
Consider the set of extensions of a rank n-1 vector bundle 
$\bar E$ by a line bundle $L$ so that 

$$deg (\bar E)-(n-1)deg (L)=s$$
The main point is to show that there is such an extension $E$
that is stable(cf[L], section 4).
We do this for a reducible curve made essentially of g elliptic curves
forming a chain. We then show that this suffices in order for the result
to hold for the generic curve. 
 
 The result was  known already when 
$rk E=2$ (cf [L-N])and when
$$0<\mu ( E'')-\mu (E')\le (g-1)/max(rkE',rk E'')$$
and in a few other special cases ([B-B-R]).

 We plan to deal with 
the general case (arbitrary rank for $E'$) in the future.

My interest on this question was spurred by a talk given by Barbara Russo
at the Europroj 96 meeting. I would like to thank her and
the organisers of the conference especially Peter Newstead for giving
me the opportunity to attend.
I am a member of the Europroj group Vector Bundles on Curves.

\bigskip
{\bf 1 The result for the special curve}
\bigskip

In this section $C$ will denote a reducible nodal curve
that we construct as follows:
Take g generic elliptic curves $C_1,...,C_g$ with marked points
$P_i, i=2,...,g; Q_i, i=1,...,g-1, P_i,Q_i \in C_i$.
Glue $C_i$ to $C_{i+1}$ at $Q_i$ and $P_{i+1}$ respectively.
We want to prove the equivalent of Theorem (0,1) for this curve.

We shall assume that the reader is familiar with the theory of
moduli of stable vector bundles on a reducible curve as developped
by Seshadri (cf.[S])(see also [T1,2,3]).

We define $\bar U_s(n,d)$ as the subset of $U_s(n,d)$ of those
$E$ such that there is a sequence (0,1) in which the $E',E''$
that give the value $s$ are both stable.
\bigskip
\proclaim (1,1)Theorem. On the curve C defined above, 
$dim\bar U_{s,1}(n,d)=d_{s,1}, 0< s\le(n-1)g, s\equiv d(n)$.

\bigskip

We introduce the following notation that we shall use repeatedly.
If $E$ is a vector bundle, $E_P$ will denote the vector space
fiber of $E$ at the point $P$ and ${\bf P}(E_P)$  the 
projective space of lines in this vector space. Also  ${\cal S}^k(E)$
denotes the set of immersed line subbundles of
degree $k$ in $E$.  

Construct now a vector bundle $E$ on $C$ of degree $d, 0\le d<n-1$
in the following way:
Denote by $k$ the greatest common divisor of $d$ and $n$. Write
$d=kd_1, n=kn_1$. On the component $C_1$ choose a vector bundle
$E_1$ obtained as the direct sum of $k$ irreducible generic 
vector bundles of degree $d_1$ and rank $n_1$. On the components
$C_i,i\ge 2$, choose a direct sum of $n$ generic line bundles
of degree zero.
Choose a positive integer $a, 0\le a\le g-1$.Denote by $[x]$ the greatest
integer less than or equal to $x$.  Define integers $t,k$
as follows if $a\le g-1-n+d$

$$(1.2.a)  t=[ (g-1-a-n+d)/n], k=g-1-a-n+d-tn $$
and if $a>g-1-n+d$
$$(1.2.b), t=0, k=g-1-a$$

The gluings will be generic on the first $g-1-a$ nodes. On the last $a$ 
nodes they will be as described below.

Pick a degree zero line summand of $E_i, i=g-a+1...g$ and glue them 
together to make a degree zero line subbundle in the chain consisting 
of the last $a$ curves. Let the gluings be generic otherwise.
 We shall see the following

\bigskip
{\bf (1,3) Claim} If $a\le g-1-n+d$, there are line subbundles of degree  
$-(k+(n-1)t+n-d-1)$ on the union of the first $g-a$ components.
 Their fibers at $Q_{g-a}$ vary in a subvariety of dimension $k$ 
of ${\bf P}(E_{g-a,Q_{g-a}})$.
If $a>g-1-n+d$, there are line subbundles of degree  
$-k$ on the union of the first $g-a$ components.
 Their fibers at $Q_{g-a}$ vary in a subvariety of dimension $k+d$ 
of ${\bf P}(E_{g-a,Q_{g-a}})$.
\bigskip

Glue one of the special directions whose existence is claimed in (1,3)
with the direction of the degree zero line subbundle. Let the gluing 
be generic otherwise.
\bigskip

{\bf (1,4) Claim.} 
$$s(E)=d+n(n-d+(n-1)t+k-1), a\le g-1-n+d$$
$$s(E)=d+nk, a>g-1-n+d$$

\bigskip
{\bf Remark}: The set of s obtained in that way cover the whole range 
of integers $s\equiv d(n), 0\le s\le (n-1)g$.
\bigskip

Let us now count the number of moduli of such families. 
The choice of $E_i$ is the most generic inside one component of
the moduli space of vector bundles on $C$(cf [T]p.342, Theorem).
The choice of gluings at each of the last $a-1$ nodes imposes 
$n-1$ conditions on the gluing. The gluing at the $g-a^{th}$
node imposes $n-1-k$ conditions if $a\le g-1-n+d$ and
$n-1-k-d$ conditions if $a> g-1-n+d$.
Therefore the number of moduli for such a family is 
$n^2(g-1)+1-(n-1)(a-1)-(n-1-k)$ if $a\le g-1-n+d$
and $n^2(g-1)+1-(n-1)(a-1)-(n-1-k-d)$  if $a> g-1-n+d$.
Using the definition (1,2) of $k$, one can check that this number 
coincides with $d_s$.

In order to prove (1,3) and finish the proof of (1,4), we need 
some preliminary results.

\bigskip

\proclaim (1,5)Lemma. Consider generic elliptic curves $C_1,...,C_t$ with
marked points $P_i,Q_i\in C_i$. Form a nodal curve C of genus t by gluing
$C_i$ to  $C_{i+1}$ at $Q_i$ and $P_{i+1}$ respectively
. Consider a vector bundle $E$
on $C$ such that the restriction $E_i$ to $C_i$ is a direct sum of $n$
generic line bundles of degree zero and the gluings
$$\varphi _i:E_{i,{Q_i}}\rightarrow E_{i+1,P_{i+1}}$$ 
at each node are generic. Consider  one-dimensional spaces $V^0_1$ of
 $E_{1,{P_1}}$ and  $W^0_t$ of $E_{t,{Q_t}}$. There is a subline bundle
 $L$ of $E$ of degree at least $-t$ such that 
 $L_{1,{P_1}}=V^0_1$, $L_{t,{Q_t}}=W_t$ if and only if $t\ge n-1$.
\bigskip

Proof. In order to construct  a line subbundle of $E$, one should take
a line subbundle $L_i$ of each $E_i$ such that the fibers at the nodes glue 
by means of the $\varphi _i$. 

As $E_i$ is a direct sum of n generic line bundles of degree zero, the
maximum degree of a line subbundle is zero and there are only a finite number
of them. As the gluings are generic, these linesubbundles do not glue
with each other.

Recall that $V^0_1$ is a given one-dimensional subspace of the vector space
$E_1$. Define the following sets:

$$A_j=\{ (L_1...L_j,V_2...V_j,W_1...W_j)\in  {\cal S}^{-1}(E_1)\times...
\times {\cal S}^{-1}(E_j)\times{\bf P}(E_{2,P_{2}})\times..
.\times{\bf P}(E_{j,P_{j}}) \times$$
$$\times{\bf P}(E_{1,Q_{1}})\times...\times{\bf P}(E_{j,Q_{j}})
|L_{1,P_1}=V^0_1,L_{i,P_i}=V_i ,i=2..j, L_{i,Q_i}=W_i,i=1...j, 
\varphi _i(W_i)=V_{i+1},i=1..j-1\}$$
\bigskip

{\bf(1,6) Claim.} $dim A_j=j,1\le j\le n-1$.
 
We prove the claim by induction on j.

For every line bundle $L_1$ of degree $-1$ on $C_1$, 
$h^0(E_1\otimes L_1^{-1})=n$. Hence, there is an (n-1)-dimensional
projective space of non-zero maps $L_1\rightarrow E_1$. If we fix one
direction  $V^0_1\subset E_{1,P_1}$, we obtain a unique such map
up to homothety. Therefore, $A_1$ is an elliptic curve isomorphic to 
$C_1$.

Assume now the result for $j-1$ and prove it for $j$.
There is a natural forgetful map $A_j\rightarrow A_{j-1}$. 
The fiber over a generic point  $(L_1...L_{j-1},V_2...V_{j-1},W_1...W_{j-1})$
is given by
 $$\{(L_j,V_j,W_j)|V_j=L_{j,P_j}=\varphi _{j-1}(W_{j-1}), 
W_j=L_{j,Q_j}\}\cong \{ L_j| \varphi _{j-1}(W_{j-1})=L_{j,P_j} \}.$$
As in the case of $A_1$, one can see that this fiber is an elliptic 
curve isomorphic to $C_j$. Therefore $dim A_j=dim A_{j-1}+1=j$. This
proves (1,6).

\bigskip
  {\bf (1,7) Claim.} Denote by $\pi _i$ projection from $A_i$ on
 the $i^{th}$ term. Then $dim\pi_{3j-1}(A_j)=j$.

Again, we prove the claim by induction on j.
For j=1, we need to see that if $L_1,L'_1$ are different generic vector bundles
of degree $-1$ on $C_1$ immersed in $E_1$ so that $L_{1,P_1}=V^0_1=L'_{1,P_1}$,
then $L_{1,Q_1}\not= L'_{1,Q_1}$.
In fact, as $Q_1$ is generic, the opposite would mean that $L_{1,Q}=L'_{1,Q}$
for all points $Q\in C_1$. Hence the two line subbundles would coincide
as bundles immersed in $E_1$. As they have the same degree, they should
in fact be equal.

Assume now the result true up to $j-1$ and let us prove it for $j$.
As $dim(A_j)=j$, if $dim(\pi_{3j-1})^{-1}(A_j)<j$ then the fibers of $\pi_{3j-1}$
are positive dimensional.

Fix $W_j\subset E_{Q_j}$. For each $L_j$, there is one immersion in
$E_j$ such that $L_{j,Q_j}=W_j$. If $dim\pi_{3j-1}(A_j)\le j-2$, 
then  $dim(\pi_j\times \pi_{3j-1})(A_j)\le j-1=dim\pi_{3(j-1)-1}(A_{j-1})=
dim  \pi_{2j-1}(A_j)$ where the last equality comes from the fact that
$W_{j-1}$ determines $V_j$ uniquely . This would imply that for a generic 
point $V_j\in {\bf P}(E_{j,P_j})$, there is only a finite number 
of line bundles $L_j$ with $L_{j,P_j}=V_j$ and this is false. Hence, if the 
fibers are positive dimensional, they are in fact one-dimensional.
Choose an irreducible component $K$ of the fiber. Then $\pi_j$ projects
it onto $Pic^{-1}(C_j)\cong C_j$. Also, at least one (say $\pi_k$)
 of the projections $\pi_i, i=1...j-1$, is not constant .
 This comes from the fact that  $L_i$ determines both $V_i$ and $W_i$
for each i.
Denote still with $\pi_k, \pi_j$ the restrictions of the projections
to $K$. Choose any base point in $C_j$. Then addition of points
of $C_j$ is well defined. We define now a non-constant map

$$\matrix {C_k&\rightarrow &C_j\cr
P&\rightarrow &\sum_{Q\in \pi_k^{-1}(P)}Q\cr}$$

Hence, $C_k,C_j$ are isogenous. This contradicts the generic choice of the
curves and finishes the proof of (1.4).

With the notations introduced so far, Lemma (1,5) amounts to saying that
 $\pi_{3j-1}(A_j)$ is a dominant map if and only if $j\ge n-1$. This
follows immediately from (1,7).
\bigskip

\proclaim (1,8) Lemma. Let all the data be as in (1.4) except that now 
$E_1$ is  the direct sum of $k$ irreducible generic 
vector bundles of degree $d_1$ and rank $n_1$. There is a line subbundle
$L$ contained in $E$  with $degL\ge -(t-1)$ such that $L_{Q_t}=W_t^0$
if and only if $t\ge n-d$.
\bigskip

Proof of (1.8):A sublinebundle of $E_1$ has degree at most 0. 
Moreover, given any line bundle $L_1\in {\bf Pic} ^0(C_1)$,
$h^0(E_1\otimes L_1^{-1})=d$. Therefore there is a projective space 
${\bf P}^{d-1}$ of immersions $L_1\rightarrow E_1$.

Consider the set 
$$\bar A_j=\{ (L_1...L_j,V_2...V_j,W_1...W_j)\in  {\cal S}^0(E_1)
\times {\cal S}^{-1}(E_2)\times...
\times {\cal S}^{-1}(E_j)\times{\bf P}(E_{2,P_{2}})\times..
.\times{\bf P}(E_{j,P_{j}}) \times$$
$$\times{\bf P}(E_{1,Q_{1}})\times...\times{\bf P}(E_{j,Q_{j}})
|L_{i,P_i}=V_i ,i=2..j, L_{i,Q_i}=W_i,i=1...j, 
\varphi _i(W_i)=V_{i+1},i=1..j-1\}$$

{\bf (1,9) Claim}. The dimension $dim \bar A_1=d=dim \pi _2 (\bar A_1)$.

Proof of the claim: There is a natural map
 $A_1 \rightarrow {\bf Pic}^0 (C_1)$ that is surjective, the fibers
being projective spaces of dimension $d-1$. From this, one can 
compute the dimension of $A_1$.

We need to show now that the generic projection has finite 
fibers. To this end we shall need the following result
\bigskip

{\bf (1,10) Claim.} We can specialise $E_1$ to a direct sum of d 
line bundles of degree one and n-d line bundles of degree zero.

Proof of (1,10). It is enough to see that an irreducible line bundle
of degree $d_1$ and rank $n_1$ can be specialised to a direct 
sum of a line bundle $L$ of degree zero and an irreducible vector
bundle $F$ of rank $n_1-1$ and degree $d_1$. 
The family of extensions of $F$ by $L$ would do the job.
\bigskip

In order to prove (1,8), it is enough to see that the fibers 
of the map to ${\bf Pic}^0(C_1)$ map to different subspaces 
by $\pi _2$. In the special case 
$$E_1=L^1\oplus ...\oplus L^d\oplus L^{d+1}\oplus...\oplus L^n,
degL^i=1,i=1...d;degL^i=0 ,i=d+1...n$$
take the fiber over $L^{d+1}$ . With the canonical coordinates in
 ${\bf P}(E_{1,Q_1})$
given by the decomposition of $E_1$, the projection of this fiber
 contains the open subset $x_{d+1}\not= 0, x_i=0, i=d+2,...n$. 
On the other hand, the projection of the fiber over
a generic $ L\not= L^i, i=1...n$ is disjoint with this set.
This proves (1,9).

Proof of (1,8): Define
 $${\bar B}_j=\{ (L_2...L_j,V_2...V_j,W_1...W_j)\in 
 {\cal S}^{-1}(E_2)\times...
\times {\cal S}^{-1}(E_j)\times{\bf P}(E_{2,P_{2}})\times..
.\times{\bf P}(E_{j,P_{j}}) \times$$
$$\times{\bf P}(E_{2,Q_{2}})\times...\times{\bf P}(E_{j,Q_{j}})
|L_{i,P_i}=V_i , L_{i,Q_i}=W_i,i=2...j, 
\varphi _i(W_i)=V_{i+1},i=2..j-1, W_j=W^0_j \}$$
 where $W^0_j$ is a fixed subspace.
By (1,6) (applied to the chain in the reverse order),
 $ dim \pi_j ({\bar B}_j)=j-1$. In order for a line subbundle as in (1,7)
to exist, we need $dim \pi_j({\bar B}_j)+dim \pi_2({\bar A}_1)\ge
dim {\bf P}(E_{1,P_1})=n-1$. Hence, $j\ge n-1$ as claimed.

\bigskip
\proclaim (1,11) Corollary.  $dim{\bar A}_j=d+j-1=
dim \pi_{3j-1}({\bar A}_j), 1\le j\le n-d$.
\bigskip

Proof. We already know the result for $j=1$. The computation of the 
dimension of ${\bar A}_j$ is identical to the proof of (1,5). Hence,
$dim \pi_{3j-1}({\bar A}_j)\le d+j-1$. But if the inequality were 
strict, as this dimension can only increase in one unit with each
$j$, then $dim \pi_{3(n-d)-1}({\bar A}{n-d})<n-1$. Then, 
$\pi _{3(n-d)-1}$ would not be surjective. On the other hand, (1.8)
implies the surjectivity of this map.
This contradiction finishes the proof of (1,11).

\bigskip
We can now finish the proof of (1,3) and (1,4). 
If $a>g-1-(n-d)$, (1,2) and (1,3) are an immediate consequence of 
(1,7) and (1,10).

Assume $a\le g-1-n+d$. We want to see that there is a line subbundle
whose degree in the chain of components is as follows

$$0,\overbrace{-1,\cdots,-1}^{n-d-1},\overbrace{0,\overbrace{-1,\cdots,-1}
^{n-1},\cdots ,0,\overbrace{-1,\cdots,-1}^{n-1}}^t,
0,\overbrace{-1,\cdots,-1}^k,\overbrace{0,\cdots,0}^a$$

By (1,8), there are such line bundles on the first $n-d$ components
and their fibers at $Q_{n-d}$ fill the whole space. One can then glue
one of them to one of the degree zero summands of $ E_{n-d+1}$.
 Using (1,5), one can 
continue this line bundle to the next $nt$ components. Finally by 
(1,11), one can continue it to the next $k$ components and the fiber
at $Q_{g-a}$ moves in the stated dimension.

For a higher degree line bundle to exist, one should replace some of the 
degree -1 line subbundles of some $E_i$ by degree 0 line subbundles. 
From (1,5) and (1,8), this is impossible.  
In order to complete the proof of (1,1), it remains to show that all bundles that appear are (semi)stable for suitably chosen polarisations.

\bigskip 

\proclaim (1,12) Lemma. Consider the vectorbundle $E$, the sublinebundle
$L$ and the quotient $E/L$ defined above. Up to tensoring the above data
with line bundles of suitable degrees, one can choose polarisations
 that make them 
semistable and even stable unless $d|n, s=0$.

Proof: Notice that the restriction of $E$ (and of course of $L$)to 
any component is semistable. We show now that the same statement 
holds for $E/L$. 

On the components $C_2,...,C_g$, $E_{|C_i}=L_1\oplus ...\oplus L_n$
where each $L_k$ has degree zero and they are generic. When $L_{|C_i}$
is one of the $L_k$ with the natural immersion, the result is obvious.
 When $L_{|C_i}$ is a degree $-1$ line bundle, its image by the immersion 
in $E_{|C_i}$ is not contained in any proper direct sum of the degreee
zero line subbundles (as the gluing data are generic). Assume the 
quotient non-stable. Let $F$ be a subbundle of the quotient of degree
$d_F$ and rank $n_F$ that contradicts semistability (i.e. 
$d_F/n_F\ge \mu (E_{|C_i}/L_{|C_i})=1/(n-1)$). 
The pull-back of $F$ to $E$ is then a subbundle of $E$ of degree $d_F-1$
and rank $n_F+1$ that contains $L_{|C_i}$. Therefore, $\pi ^{-1}(F)$ is 
not contained in any proper direct sum of the $L_i$. This implies that 
$deg\pi^{-1}(F)<0$. So $d_F\le 0$ which contradicts the assumption.

We now prove the similar result for $C_1$. We shall write $E$ for 
$E_{|C_1}$.Again the result is obvious
for $d=0$. So let us assume $0< d< n$. We shall show that the quotient
of $E$ by a sublinebundle $L$ of degree zero is indecomposable.
Assume the opposite. Then, the quotient is the direct sum of 
indecomposable bundles $F_1\oplus ...\oplus F_k$.
The exact sequence
$$0\rightarrow L\rightarrow E \rightarrow F_1\oplus ...\oplus F_k
\rightarrow 0$$
corresponds to an extension class
 
$$\delta =(\delta_1,...,\delta_k)\in H^1(F^*\otimes L)=
H^1(F_1^*\otimes L)\oplus ...\oplus H^1(F_k^*
\otimes L)$$ 
Each $\delta_i$ corresponds to an extension 
$$0\rightarrow L \rightarrow E_i \rightarrow F_i \rightarrow 0$$
and we have an exact commutative diagram 
$$\matrix{ & &0& &0& & & & \cr
 & &\downarrow& &\downarrow & & & & \cr
0&\rightarrow & L&\rightarrow &E&\rightarrow &F_1\oplus...\oplus F_k&
\rightarrow &0\cr
 & &\downarrow\Delta& &\downarrow& &\Vert & & \cr
0&\rightarrow &\overbrace{ L\oplus\cdots\oplus L}^k&\rightarrow &
 E_1\oplus...\oplus E_k &\rightarrow & F_1\oplus...\oplus F_k&
\rightarrow &0\cr
 & &\downarrow& &\downarrow& & & & \cr
 & &\overbrace{L\oplus\cdots\oplus L}^{k-1}& =
&\overbrace{L\oplus\cdots\oplus L}^{k-1}& & & & \cr
 & &\downarrow & &\downarrow & & & & \cr
 & &0 & &0 & & & & \cr}$$
where $\Delta$ is the diagonal map.

 We first prove that the map from $E$ to any direct summand 
of any $E_i$ cannot be zero: Assume this were the case. Write 
$E_i=E'_i\oplus E''_i$ with the map from $E$ to $E'_i$ being zero.
Then, $E'_i$ injects into $L\oplus \cdots \oplus L$. In
particular, it is torsion-free of slope at most 0. Chassing the
diagram above, one finds that $L$ surjects onto $E'_i$ and therefore
$E'_i$ has rank one. It follows that $E'_i=L, E''_i=F_i$.
This contradicts the commutativity of the upper left square.

 As $E$ is semistable and $h^0(E^*\otimes E'_i)>0$,
 $\mu (E_i)-\mu (E)\ge 0$  (cf [A]). Equivalently,
$ndeg(E_i)-d rk(E_i)\ge 0$. Notice that $degE_i=degF_i, rk E_i=rkF_i+1$.
Adding over all inequalities above, we find
$nd-d(n-1+k)\ge 0$. So, $1-k\ge 0$ or $k\le 1$. Hence, $k=1$ and $F$ is 
irreducible.

Consider now a line bundle $L$ on $C$ of degree $k_i$ on the component
$C_i$. When tensoring  a vector bundle $E$ by $L$, the degree of $E$ is
 modified
in $nk_i$. Therefore, for a suitable choice of the $k_i$, one can bring 
the degree of the new vector bundle on each component
inside any given interval of length
n. If the restriction of $E$ to each component is semistable, so is the 
restriction of the modified vector bundle. By [T1] Theorem (see 
also [T2]), this is enough to insure stability. 
\bigskip

{\bf The reduction step}
\bigskip

\proclaim (2,1)Lemma. Let $C_0$ be (a possibly reducible) curve such that 
$\bar U_{s,n'}(n,d)(C_0)$ is non-empty and has a component of dimension 
$d_{s,n'}$. Denote by $C_1$ the curve obtained from $C_0$ by adding chains 
of rational components between the nodes. Then $\bar U_{s,n'}(n,d)(C_1)$
is non-empty and has a component of dimension $d_{s,n`}$

Proof: Consider the map $\pi :C_1\rightarrow C_0$ obtained by contracting 
the rational components. Denote by $E_1$ the pull-back by $\pi$ of a vector
 bundle $E_0$ on $C_0$. This is trivial (i.e. of the form ${\cal O }^n$)
on the additional chains of components. This establishes 
a one to one correspondence between vector bundles
on $C_0$ and vector bundles on $C_1$ that are trivial on the 
additional chains of rational components (cf[T3] section 4). This correspondence
preserves stability by suitable polarisations (cf [T3] Lemma (1.3)).
The subbundles of a trivial vector bundle have degree at most 0.
Also the group of automorphism of a trivial vector bundle is 
the linear group. This establishes an isomorphism between any two
fibers. If the chain  glues to the rest at points $Q_i$, $P_{i+1}$
coming from $C_0$, this isomorphism coincides with the original gluing
on $C_0$.One can then see that $s_{n'}(E_0)=s_{n'}(E_1)$.

Moreover in a neighborhood of $E_1$, the elements of the moduli
space are trivial on the rational chain (as these are the most generic 
elements in the moduli space). Hence the dimension of $\bar U_{s,n'}(n,d)(C_1)$
in a neighborhood of $E_1$ coincides with the dimension of 
$\bar U_{s,n'}(n,d)(C_0)$ in a neighborhood of $C_0$. This completes the proof 
of (2,1).

\proclaim (2,2) Proposition. Let $\cal C\rightarrow S$ be a family of curves
parametrised by a discrete valuation ring. Assume that the generic 
curve $C_t$  is irreducible non-singular  while the special curve 
$C_0$ is a curve of compact type (i.e. tree-like).
Assume that $C_0$ has a vector bundle $E_0$ of rank $n$ and degree $d$
that can be written as an extension 
$$0\rightarrow E'_0\rightarrow E_o\rightarrow E''_0\rightarrow 0$$
Assume that (up to tensoring with line bundles on $C_0$ of preassigned degree
on each component) all data are stable (by  polarisations
 $(a'_i),(a_i), (a''_i)$ on 
$C_0$ if the curve is reducible).Assume that $s_{n'}(E)=s$
and $E'$ is a subbundle of rank $n'$ with maximal degree.
Assume moreover than in a neighborhood of $E_0$ The locus of vector bundles
with  $s_{n'}(E)=s$ has dimension $d_{s,n'}$. Then in the generic curve, the 
locus of stable bundles with $s_1(E)=s$ is non-empty and has dimension 
$d_{s,n'}$.
\bigskip

Proof: Make a base change if necessary so that $E'_0,E_0,E''_0$ can be 
extended to vector bundles on the whole family. One can also assume
 that the family has enough sections
so that moduli spaces of stable bundles on the family exist.
 One may need to replace
the central fiber by a new curve $C_1$ which has strings of rational curves
inserted between the nodes. By (2,1), the hypothesis on $C_0$
hold also for $C_1$.

Consider now divisors on $\cal C$ consisting of sums of components of
$C_1$ with suitable multiplicity. As $C_1$ is of compact type, 
there are such divisors with any preassigned degree on the components
of $C_1$. Choose line bundles ${\cal L}',{\cal L},{\cal L}''$,
of this form on $\cal C$.
Consider now the moduli spaces $U(n',d',(a'_i))$,
$U(n'',d'',(a''_i))$. The elements of these moduli spaces tensored 
with ${\cal L}'$ and $ {\cal L}''$ respectively parametrise
vector bundles of the type of $E'$ and $E''$ that we started with.
Then, as in [L] section 4, one can consider families of extensions
of this type of bundles. 
One sees then that the family of stable extensions as stated is 
either empty or has dimension $d_s+dimS$. As the fiber over $C_1$ is 
non-empty and has dimension $d_{n',s}$, the family of stable extensions 
projects onto $S$ and the generic fiber has dimension $d_{n',s}$.

\bigskip 
Proof of (0,3): this is a consequence of (1,1) and (2,2).
\bigskip
{\bf References}
\bigskip
[A] M.Atiyah, Vector bundles on elliptic curves, Proc.London Math.
Soc. {\bf(3)7}(1957), 414-452.

[B-B-R] E.Ballico, L.Brambila-Paz, B.Russo, Exact sequences of stable
vector bundles on projective curves. Math.Nach. To appear.

[L]H.Lange, Zur Klassikation von Regelmannigfaltigkeiten,
 Math.Ann. 
{\bf 262} (1983), 447-459.

[L-N] H.Lange, M.S.Narasimhan, Maximal subbundles of rank two vector
bundles on curves. Math.Ann. {\bf 266} (1983), 55-72.

[M-S] S.Mukai, F.Sakai, Maximal subbundles of vector bundles on a curve,
Man.Math. {\bf 52} (1985), 251-256.

[N]M.Nagata, On self-intersection number of vector bundles of rank two on
 a Riemann surface, Nagoya Math. J. {\bf 37} (1970), 191-196.

[S] C.S.Seshadri. Fibres vectoriels sur les courbes algebriques. 
Asterisque {\bf 96},(1982).

[T1] M.Teixidor i Bigas, Moduli spaces of (semi)stable bundles on tree-
like curves, MathAnn.{\bf 290} (1991), 341-348.

[T2]--, Moduli spaces of vector bundles on reducible curves, Amer.J.of Math.
{\bf 117} (1995), 125-139.

 [T3]--, Compactifications of moduli spaces of vector bundles on
singular curves: two points of view. Preprint.
\bigskip 

montserrat teixidor i bigas,

 Tufts University, Medford, MA 02155, U.S.A.

Temporary address DPMMS, 16 Mills Lane, Cambridge CB2 1SB, England

teixidor@dpmms.cam.ac.uk   or mteixido@pearl.tufts.edu
\end